%
\documentstyle{amsppt} 
\magnification=\magstep1

\NoRunningHeads
\NoBlackBoxes
\parindent=1em 
\vsize=7.4in



%
\topmatter
\title
every frame is a sum of three (but not two) orthonormal bases - and other frame representations                
\endtitle
\author
Peter G. Casazza
\endauthor
\address
Department of Mathematics,
The University of Missouri,
Columbia, Missouri 65211
\endaddress
\email
pete\@casazza.math.missouri.edu
\endemail
\thanks
This research was supported by NSF DMS 9701234.  Part of this research was 
conducted while the author was a visitor at the ``Workshop on Linear Analysis
and Probability'', Texas A\&M University.
\endthanks


\abstract
We show that every frame for a Hilbert space $H$ can be written as a (multiple
of a) sum of three orthonormal bases for $H$.  We next show that this 
result is best possible by including a result of N.J. Kalton: A frame
can be represented as a linear combination of two orthonormal bases
if and only if it is a Riesz basis.  We further show that
every frame can be written as a (multiple of a) sum of two tight frames with
frame bounds one or a sum of an orthonormal basis and a Riesz basis for $H$.
Finally, every frame can be written as a (multiple of a) average of
two orthonormal bases for a larger Hilbert space. 
\endabstract
\endtopmatter
\document
\baselineskip=15pt
\heading{1.Frames as Operators}
\endheading
\vskip10pt

If $H$ is a Hilbert space, we denote the set of all bounded operators 
$T:H\rightarrow H$ by $B(H)$.  We will always use $(e_{n})$ to denote
an orthonormal basis on $H$.  Recall that a sequence $(x_{n})$ in a Hilbert space $H$ is 
called a {\bf frame} for $H$ if there are constants 
$0 < A \le B$ so that for all $x\in H$ we have
$$
A\|x\|^{2} \le \sum_{n}|<x,x_{n}>|^{2} \le B\|x\|^{2}.
$$
We call $A,\ B$ the {\bf frame bounds} for the frame and if $A = B$, we call 
this a {\bf tight frame}.  The frame definition has many equivalent forms.  We
will work here with frames thought of as operators on $H$.  That is, a sequence 
$(x_{n})$ is a frame on $H$ if and only if there is an operator $T:H\rightarrow H$ so that
$Te_{n} = x_{n}$ and $T$
is an onto map.  This equivalence is easily checked.  For one direction,
given the operator $T$ we just check that 
$$
T^{*}(x) = \sum_{n}<x,x_{n}>e_{n}.
$$
Since $T$ is onto, $T^{*}$ is an (into) isomorphism.  Hence, $TT^{*}$ is an
onto isomorphism and 
$$
TT^{*}(x) = \sum_{n}<x,x_{n}>x_{n},
$$
is the so called {\bf frame operator}.  Conversely, if $(x_{n})$ is a frame,
then $(x_{n})$ is a Hilbertian sequence in $H$.  That is, for all sequences
of scalars $(a_{n})$ we have
$$
\sum_{n}a_{n}x_{n}
$$
converges in $H$.  Hence, if we define $T:H\rightarrow H$ by $Te_{n} = x_{n}$,
then $T$ is a bounded linear operator from $H$ to $H$.  Also, $TT^{*}$ is 
the frame operator for this frame and hence $TT^{*}$ is an onto isomorphism.
Therefore, $T$ must be an onto map.  

So we can consider the ``equivalence'' between frames and onto maps on $H$.  
That is, if we have a theorem about onto maps on $H$ (or about all bounded
operators on $H$), then we have a theorem about frames.  

\heading{2.  Frames as Sums}
\endheading
\vskip10pt

All Hilbert spaces will be complex (at the end we will discuss what
happens in the real case).  The results we
use from Operator Theory can be found in any standard book in Functional
Analysis \cite{2}, or books on Hilbert space theory \cite{1,3}.  
Recall that a {\bf unitary operator} $U:H\rightarrow H$ is an onto isometry, a {\bf partial isometry} is 
an operator which is an isometry on the orthogonal complement of its kernel, a
{\bf co-isometry} is an operator whose adjoint is an into isometry, and
a {\bf maximal partial isometry} is either an isometry or a co-isometry.  An
operator $T\in B(H)$ is called a 
{\bf positive operator} if for all $x\in H$ we have $<Tx,x>\ge 0$.  That is,
$<Tx,x>$ is both a real number and positive. The main result we will use is the fact that every $T\in B(H)$ has a representation in the form $T = VP$ (called 
the{\bf Polar decomposition of T}) where $V$ is a maximal partial isometry and
$P$ is a positive operator.  Moreover, we may assume that ker $V$ = ker $P$.  
Also, every positive operator $P$ on $H$ with $\|P\|\le 1$ can be written
in the form $P = \frac{1}{2}(W+W^{*})$, where $W = P+i\sqrt{1-P^{2}}$ is unitary.   
  
Our first result is well known in operator theory, but since it is never 
formally stated, we include a sketch of the proof. 
\proclaim{Proposition 2.1}
If $T:H\rightarrow H$ is a bounded linear operator, then $T = a(U_{1}+U_{2}+U_{3})$, where each $U_{j}$ is a unitary operator and for any $0 < \epsilon < 1$, we can specify $a = \frac{\|T\|}{1-\epsilon}$.
\endproclaim

\demo{Proof}
Fix $0 < \epsilon < 1$ and let
$$
S = \frac{1}{2}I+\frac{1-\epsilon}{2}\frac{T}{\|T\|}.
$$
Then
$$
\|I-S\| = \|\frac{1}{2}I-\frac{1-\epsilon}{2}\frac{T}{\|T\|}\| \le \frac{1}{2}+\frac{1-\epsilon}{2}< 1.
$$
It follows that $S$ is an onto isomorphism.  We now write the polar decomposition of $S$ as $S = VP$ where $V$ is a maximal partial isometry and $P$ is a positive operator and $\text{ker}\ V = \text{ker}\ P$.  Since $S$ is an isomorphism, $P$ is an isomorphism and $V$ is a unitary operator (Since a necessary and sufficient condition that $V$ be an isometry is that $S$ be 1-1, and a necessary and sufficient condition that $V$ be a co-isometry is that S has dense range \cite{3}).
  Also, $\|S\| < 1$ implies that $\|P\|\le 1$, and hence $P=\frac{1}{2}(W+W^{*})$.
Now we have the representation,
$$
S = \frac{1}{2}(VW+VW^{*}),
$$
where $VW,\ VW^{*}$ are unitary.  Finally,
$$
T = \frac{\|T\|}{1-\epsilon}(VW+VW^{*}-I)
$$
is the required decomposition of $T$.
\enddemo

As an immediate consequence of Proposition 2.1 we have,

\proclaim{Corollary 2.2}
If $(x_{i})_{i\in I}$ is a frame for a Hilbert space $H$ with upper frame bound $B$, then for every ${\epsilon}> 0$, there are orthonormal
bases $(f_{i}),\ (g_{i}),\ (h_{i})$ for $H$ and a constant $a=B(1+\epsilon)$ so that
$$
x_{i} = a(f_{i}+g_{i}+h_{i}),\ \ \forall \ i\in I.
$$
\endproclaim

\demo{Proof}
If we choose an orthonormal basis $(e_{i})$ for $H$, then
the operator $T:H\rightarrow H$ given by $T(e_{i}) = x_{i}$ is an onto map and
$\|T\| = B$.  By Proposition 2.1, we can write $T = a(U_{1}+U_{2}+U_{3})$ where
each $U_{j}$, $1\le j\le 3$ is a unitary operator.  It follows that $(U_{j}e_{i})_{i\in I}$ is an orthonormal basis for $H$ for $1\le j\le 3$ and by the proof of Proposition 2.1 (with $\frac{1}{1-\epsilon}$ there equaling 
$1+\epsilon$ here)we get the estimate for $a$.
\enddemo

Corollary 2.2 cannot be improved to represent the frame as a sum of only
two orthonormal bases (or even orthonormal sequences) in general, as the next example shows. We call a frame a {\bf normalized tight frame} if it is a tight frame with frame bounds 1.  Even though the proposition following the
example supersedes it, we have included it here since it works in both
the real and complex case and clearly illuminates the difficulties
encountered in such a representation.  

\proclaim{Example 2.3}
There is a normalized tight frame for a (real or complex) Hilbert space
$H$ which cannot be written 
as any linear combination
of two orthonormal sequences in $H$.
\endproclaim

\demo{Proof}
Let $(e_{i})$ be an orthonormal basis for $H$ and consider the normalized
tight frame: $x_{1} = 0,$ and for all $1\le i$, $x_{i+1} = e_{i}$.  We proceed by way of contradiction.  If we can find orthonormal sequences $(f_{i}),\ (g_{i})$ in $H$ and numbers $a,b$ so that $x_{i} = af_{i}+bg_{i}$, for all $i\in I$ then
$$
x_{1} = 0 = af_{1}+bg_{1}.
$$
Hence, if $a\not= 0\not= b$, then $\text{span}\ f_{1} = \text{span}\ g_{1}$ and orthogonality imply
$$
\text{span}(f_{i})_{i=2}^{\infty} = \text{span}(g_{i})_{i=2}^{\infty}\not= H,
$$
while
$$
\text{span}(af_{i}+bg_{i})_{i=2}^{\infty} = \text{span}(e_{i})_{i=1}^{\infty} = H.
$$
This contradiction completes the proof of the Example.  
\enddemo

Note that the above example does not even allow us to find two sequences in
$H$ which are only orthonormal bases for their spans and add up to our frame.

The following result, due to N.J. Kalton and appearing here with his
permission, gives a complete characterization of those frames which
can be written as linear combinations of two orthonormal bases as
precisely the class of Riesz bases.  Again, it comes from a result in
operator theory.

\proclaim{Proposition 2.4}
If $T\in B(H)$ is onto, then $T$ can be written as a linear combination
of two unitary operators if and only if $T$ is invertible.
\endproclaim

\demo{Proof}
$\Leftarrow$:  If $T$ is invertible, then by the proof of Proposition 2.1,
we have that $T = a(U_{1}+U_{2})$, where $U_{1},U_{2}$ are unitary operators.

$\Rightarrow$:  Suppose $T = aU_{1}+bU_{2}$, where 
$U_{1},U_{2}$ are unitary operators.  If either of $a$ or $b$ equals zero,
we are done.  So without loss of generality, (after dividing by the smaller
of the two and observing that $T$ is invertible if and only if $\frac{1}{a}T$
is invertible) we may assume that
$$
T = U_{1} + aU_{2},
$$
where $|a|\ge 1$.  For all $0\le t <1/2$, let
$$
S_{t} = tU_{1} + (1-t)aU_{2}.
$$
We observe that for all $0\le t < 1/2$, the operator $S_{t}$ is a 
(possibly into) isomorphism.  To see this we calculate for all $f\in H$,
$$
\|S_{t}(f)\| \ge (1-t)|a|\|U_{2}(f)\|-t\|U_{1}(f)\| = 
[(1-t)|a| -t]\|f\|.
$$
Since $0\le t<1/2$ and $|a|\ge 1$, we have that $[(1-t)|a| -t]>0$ and
so $S_{t}$ is an isomorphism.  Also, by our assumption that $T$ is onto,
we have that $S_{1/2}$ is onto.  A result of Kalton, Peck, and
Roberts \cite{5}, Proposition 7.8 (this is done for open maps but works perfectly well for onto maps) we have that the onto maps for an open set 
in $B(H)$.  Since 
$$
\lim_{t\rightarrow 1/2}S_{t} = \frac{1}{2}T,
$$
it follows that for $t$ close enough to $1/2$, the operator $S_{t}$ is onto.
Hence $S_{t}$ is invertible (being an isomorphism).  But Proposition 7.9 of
\cite{5} says that the invertible operators forms a clopen (both closed and
open) subset of the onto maps in $B(H)$.  Therefore, $S_{1/2}$ is also an
invertible operator.  Hence, $T = 2S_{1/2}$ is invertible.  

\enddemo

The corresponding result for frames looks like:

\proclaim{Proposition 2.5}
A frame $(f_{i})$ for $H$ can be written as a linear combination of
two orthonormal bases for $H$ if and only if $(f_{i})$ is a Riesz basis
for $H$.
\endproclaim

Although an arbitrary frame cannot be written as a multiple of a sum of two
orthonormal sequences in the space, it can be written as a multiple of a 
sum of two normalized tight frames in the space. This follows from another
general result from operator theory concerning 
the decomposition of an operator.
That is, every operator $T$ on a Hilbert space can be written in the form
$$
T = VP = \frac{\|T\|V}{2}(W+W^{*}), \tag 2.1
$$
where $W$ is unitary and $V$ is a maximal partial isometry.  It follows that 
$VW$ and $VW^{*}$ are maximal partial isometries.  That is, each of these operators is either an isometry or a co-isometry.  However, if $T$ induces
a frame on $H$ then $T$ has dense range and so $V$ must be a co-isometry.    

\proclaim{Proposition 2.6}
If $T$ is a co-isometry on $H$, and if $(e_{i})$ is a orthonormal basis
 for $H$, then $(Te_{i})$ is a normalized tight frame for $H$.
\endproclaim

\demo{Proof}
Since $T$ is a co-isometry, $T^{*}$ is an isometry.  Hence, for all $f\in H$,
$$
\|f\|^{2} = \|T^{*}f\|^{2} = \sum_{i}|<T^{*}f,e_{i}>|^{2} 
=\sum_{i}|<f,Te_{i}>|^{2}.
$$
Therefore, $(Te_{i})$ is a normalized tight frame for $H$. 
\enddemo

From our discussion preceeding Proposition 2.6, we see that any onto map
$T:H\rightarrow H$ (inducing the frame $(Te_{i})$ on $H$) can be written in the form given in (2.1) where each of $VW,\ VW^{*}$ is 
a co-isometry.  This combined with Proposition 2.6 gives immediately

\proclaim{Proposition 2.7}
Every frame for a Hilbert space is (a multiple of) the sum of two normalized 
tight frames for $H$.
\endproclaim

Proposition 2.7 should be compared to the fact that every frame is
equivalent to a normalized tight frame.  As we saw earlier, even normalized
tight frames may not be the sum of two orthonormal sequences in $H$.  But,
if we are willing to ``expand'' the Hilbert space, then we can get a good
representation.  The main ideas below are part of a nice alternate 
approach to frames due to Han and Larson \cite{4}.

\proclaim{Proposition 2.8}
If $(x_{i})$ is a normalized tight frame for $H$, then there is a Hilbert
space $H\subset K$ and two orthonormal bases $(f_{i}),\ (g_{i})$ for $K$
so that $(x_{i})$ is the average of $(f_{i})$ and $(g_{i})$.
\endproclaim

\demo{Proof}
Let $Q:H\rightarrow H$ be the onto map $Qe_{i} = x_{i}$.  Then $(x_{i})$ is
a normalized tight frame implies that $Q$ is a quotient map and hence $Q^{*}$
is an isometry.  By identifying $x_{i}$ with $Q^{*}x_{i}$ and $H$ with $Q^{*}(H)$
we can let $K = H$ and choose the orthogonal projection $P:K=H\rightarrow Q^{*}(H)$ and observe that $P(e_{i}) = Q^{*}x_{i}$.  That is, for all $i,j$
 we have:
$$
<Q^{*}x_{i},Pe_{j}> = <PQ^{*}x_{i},e_{j}> = <Q^{*}x_{i},e_{j}> = 
<x_{i},Qe_{j}> 
$$
$$
= <x_{i},x_{j}> = <Q^{*}x_{i},Q^{*}x_{j}>.
$$
But, $(Q^{*}x_{i})$ spans $Q^{*}H$ and so $Pe_{j} = Q^{*}x_{j}$.  Finally,
$(Pe_{i}+(I-P)e_{i})$ and $(Pe_{i} - (I-P)e_{i})$ are both orthonormal bases
for $K$ and their average is $Pe_{i} = Q^{*}x_{i}$. 
\enddemo

Since every frame is equivalent to a normalized tight frame we have:

\proclaim{Corollary 2.9}
Every frame is equivalent to a frame which is an average of two orthonormal
bases for a larger Hilbert space. 
\endproclaim

Riesz bases are a little more general than orthonormal bases, and so we
get a stronger result relative to frames.

\proclaim{Proposition 2.10}
Every frame for a Hilbert space $H$ is (a multiple of) the sum of a
orthonormal basis for $H$ and a Riesz
basis for $H$. 
\endproclaim

\demo{Proof}
We proceed as in Proposition 2.1 with a slight change.  Given a Pre-frame
operator $T:H\rightarrow H$ with our frame being $(Te_{i})$, define an
operator $S$ by
$$
S = \frac{3}{4}I + \frac{1}{4}(1-\epsilon)\frac{T}{\|T\|}.
$$
Then again we have $\|I-S\| < 1$ and $\|S\|\le 1$, so $S$ is an invertible operator and
as in the proof of Proposition 2.1 we can write
$$
S = \frac{1}{2}(W+W^{*}),
$$
where $W$ is a unitary operator.  (Note that here $W$ is taking the
place of $VW$ in Proposition 2.1).  Now we have:
$$
T = \frac{4\|T\|}{(1-\epsilon)}\left [ \frac{1}{2}(W+W^{*}) - \frac{3}{4}I \right ] .
$$
Hence,
$$
T = \frac{2\|T\|}{(1-\epsilon)}[W + R], \ \ \ \text{where}\ \ 
R = W^{*}-\frac{3}{2}I.
$$
Now, $W$ is unitary so $(We_{i})$ is an orthonormal basis, and $W^{*}$ is
unitary implies that $R$ is an isomorphism (possibly into).  But, it is
easily checked that $R$ is onto since 
$$
\|I - \frac{-1}{2}R\| = \| \frac{1}{4}I + \frac{1}{2}W^{*}\| < 1.
$$
Thus $\frac{-1}{2}R$
is an invertible operator and hence $R$ is an invertible operator.  Since $R$
is an invertible operator, $(Re_{i})$ is a Riesz basis for $H$.
\enddemo

\proclaim{Remark 2.11}
Again, what we have really used in Proposition 2.10 is just a result from
operator theory which says that every operator on a Hilbert space is a
multiple of the sum of a unitary operator on $H$ and an invertible 
operator on $H$.
\endproclaim

\proclaim{Remark 2.12}
In the real case, we cannot write a positive operator as an average of two
unitaries.  This result in the complex case comes from the representation
of the extreme points of the ball of $B(H)$ and the fact that every positive operator $P$ with $\|P\| = 1$ is actually an average of two extreme points.  In the real case, we lose this representation but do have a similar one with
a representation in terms of 16 operators. Thus, we can recapture the above results with "larger sums".  For
example, we can write a frame in a real Hilbert space 
as a multiple of a sum of 16 orthonormal bases 
or 4 Riesz bases, and a normalized tight frame can be written as the sum
of four orthonormal bases etc.
\endproclaim

\enddocument